\documentclass[12pt]{article}

    \usepackage{psfrag}

    \usepackage{graphics}

    \usepackage{epsfig}
    \usepackage{amsthm,amsfonts,amsmath}
    \usepackage{color}   

    \usepackage[dvips]{geometry}
  \newcommand{\R}{\mathbb{R}}

\newtheorem{teorema}{Theorem}[section]
\newtheorem{lema}{Lemma}[section]
\newtheorem{prop}{Proposition} [section]

\newtheorem{obs}{Remark} [section]

\newtheorem{defi}{Definition} [section]

\begin {document}

\begin{center}\Large{ \textbf{On $H=1/2$ Surfaces in $\widetilde{PSL}_{2}(\R,\tau)$}}
\end{center}

\vspace{0.2cm}

\begin{center}\large{Carlos Pe\~nafiel \footnote{The author was partially supported by FAPERJ - Brasil.\\ e-mail: carnilexe@yahoo.com.br}}
\end{center}
\vspace{0.8cm}

\textbf{Abstract.} We study the behavior of $H=1/2$ surfaces immersed in $\widetilde{PSL}_{2}(\R,\tau)$.

\vspace{1cm}
\vspace{0.5cm}

\textbf{Acknowledgements.} The author would like to thank to the Professor Harold
Rosenberg, for suggesting this problem, for many interesting and
stimulating discussion on the subject and for the constant support
throughout his work.

\section{Introduction}

\indent

In this paper we study complete constant mean curvature $H=1/2$ surfaces in $\widetilde{PSL}_{2}(\R,\tau)$. Recall that in \cite{HRS} the authors generalized to $\mathbb{H}^{2}\times\R$ the half-space theorem of Hoffman and Meeks which says that a properly immersed minimal surface in $\R^3$ that lies in a half-space must be a plane. The main theorem in \cite{HRS} says that a properly embedded constant mean curvature $H=1/2$ surface in $\mathbb{H}^2\times\R$ which is asymptotic to a horocylinder $C$ and on one side of $C$; such that the mean curvature vector of the surface has the same direction as that of $C$ at points of the surface converging to $C$, then the surface is equal to $C$ (or a subset of $C$ if the surface has non-empty boundary).

We follow these ideas to show an analogous result in the space $\widetilde{PSL}_{2}(\R,\tau)$. More precisely, our main theorem says:
\begin{teorema} \label{01} Let $\Sigma$ be a properly embedded constant mean curvature $H=1/2$ surface in $\widetilde{PSL}_{2}(\R,\tau)$. Suppose $\Sigma$ is asymptotic to a horocylinder $C$, and on one side of $C$. If the mean curvature vector of $\Sigma$ has the same direction as that of $C$ at points of $\Sigma$ converging to $C$, then $\Sigma$ is equal to $C$ (or a subset of $C$ if $\partial\Sigma\neq\phi$)
\end{teorema}

This theorem lets us obtain (in the same sense as in \cite{HRS}) the following result:
\begin{teorema} Let $\Sigma$ be a complete immersed surface in $\widetilde{PSL}_{2}(\R,\tau)$ of constant mean curvature $H=1/2$. If $\Sigma$ is transverse to $E_3=\partial_t$ then $\Sigma$ is an entire vertical graph over $\mathbb{H}^2$.
 \end{teorema}

\section{The Laplacian and the Mean Curvature Equation for $H$ Surfaces in $\widetilde{PSL}_{2}(\R,\tau)$}

\indent

The space $\widetilde{PSL}_{2}(\R,\tau)$ is a complete simply connected homogeneous manifold. Such a manifold is a Riemannian fibration over the 2-dimensional Hyperbolic space $\mathbb{H}^2$. That is, there exist a Riemannian submersion $\pi:\widetilde{PSL}_{2}(\R,\tau)\rightarrow\mathbb{H}^2$, which also is a Killing submersion (see Definition \ref{a1}).

The space $\widetilde{PSL}_{2}(\R,\tau)$ is topologically $\mathbb{H}^2\times\R$ each fiber is diffeomorphic to $\R$ (the real line) and the bundle curvature of the submersion is $\tau$.

It is well know, that the space $\widetilde{PSL}_{2}(\R,\tau)$ is given by (see \cite{D})

\begin{equation*}
\widetilde{PSL}_{2}(\R,\tau)=\{(x,y,t)\in\R^3;y>0\}
\end{equation*}
endowed with the metric
\begin{equation*}
g=\lambda^2(dx^2+dy^2)+(-2\tau\lambda dx+dt)^2, \hspace{0.2cm} \lambda=\dfrac{1}{y}.
\end{equation*}

There is a natural orthonormal frame $\{E_1,E_2,E_3\}$ given by (in coordinates $\{\partial_x,\partial_y,\partial_t\}$),
\begin{equation*}
 E_1=\dfrac{\partial_x}{\lambda}+2\tau\partial_t, \hspace{0.3cm} E_2=\dfrac{\partial_y}{\lambda}, \hspace{0.3cm} E_3=\partial_t.
\end{equation*}
$E_3$ is the killing field tangent to the fibers. The metric $g$ induces a Riemannian connection $\overline{\nabla}$ given by:

\begin{center}
  $\begin{array}{ccc}
    \overline{\nabla}_{E_{1}}E_{1}=-\displaystyle\frac{\lambda_{y}}{\lambda^{2}}E_{2} \ \ \ \ \ \ \ \ \ \ \ \ & \overline{\nabla}_{E_{1}}E_{2}=\displaystyle\frac{\lambda_{y}}{\lambda^{2}}E_{1}+\tau E_{3} & \overline{\nabla}_{E_{1}}E_{3}=-\tau E_{2} \\[10pt]
    \overline{\nabla}_{E_{2}}E_{1}=\displaystyle\frac{\lambda_{x}}{\lambda^{2}}E_{2}-\tau E_{3} & \overline{\nabla}_{E_{2}}E_{2}=-\displaystyle\frac{\lambda_{x}}{\lambda^{2}}E_{1} & \overline{\nabla}_{E_{2}}E_{3}=\tau E_{1} \\[10pt]
    \overline{\nabla}_{E_{3}}E_{1}=-\tau E_{2} \ \ \ \ \ & \overline{\nabla}_{E_{3}}E_{2}=\tau E_{1} & \overline{\nabla}_{E_{3}}E_{3}=0 \\[10pt]
    [E_{1},E_{2}]=\displaystyle\frac{\lambda_{y}}{\lambda^{2}}E_{1}-\frac{\lambda_{x}}{\lambda^{2}}E_{2}+2\tau E_{3} & [E_{1},E_{3}]=0 & [E_{2},E_{3}]=0 \\[10pt]
  \end{array}$
\end{center}

Following ideas of \cite{HRS}, we will consider horizontal graphs $y=f(x,t)$. Denoting by $S=graph f$, we have the following lemma:
\begin{lema}
 Denoting by $H$ the length of the mean curvature vector of $S$. Then, the function $f$ satisfies the equation
\begin{equation*}
 2H\lambda^2W^3=(f^2+f_t^2)f_{xx}-2(f_xf_t-2\tau f)f_{xt}+((1+4\tau^2)+f_x^2)f_{tt}+f(1+f_x^2)+2\tau f_xf_t
\end{equation*}
where $W=\sqrt{f^2+f^2_t+f^2(f_x+2\tau\lambda f_t)^2}$.\\
In particular the horocylinders $f(x,t)=cte$ have constant mean curvature $H=1/2$.
\end{lema}
\begin{proof}
 The surface $S=gragh f$ is parameterized by $\varphi(x,t)=(x,f(x,t),t)$, so the adapted frame to $S$ is given by:
\begin{eqnarray*}
 \varphi_x &=& \lambda(E_1+f_x E_2-2\tau E_3) \\
 \varphi_t &=& \lambda f_tE_2+E_3 \\
 N &=& \dfrac{-(f_x+2\tau\lambda f_t)E_1+E_2-\lambda f_t E_3}{\sqrt{1+(f_x+2\tau\lambda f_t)^2+\lambda^2f_t^2}}
\end{eqnarray*}

$N$ is the unit normal to $S$. Denoting by $g_{ij}$ and $b_{ij}$ the coefficients of the first and second foundamental form respectively we have that the function $H$ satisfies the equation:
\begin{equation*}
 2H=\dfrac{b_{11}g_{22}+b_{22}g_{11}-2b_{12}g_{12}}{g_{11}g_{22}-g_{12}^2}
\end{equation*}
 Since,
\begin{eqnarray*}
 \overline{\nabla}_{\varphi_x}\varphi_x &=& -\lambda^2f_x(2+4\tau^2)E_1+[\lambda f_{xx}+\lambda^2((1+4\tau^2)-f_x)]E_2 +2\tau\lambda^2f_xE_3 \\
  \overline{\nabla}_{\varphi_t}\varphi_x &=& [\tau\lambda f_x=\lambda^2f_t(1+2\tau^2)]E_1+[\lambda f_{xt}-\lambda^2f_xf_t-\lambda\tau]E_2 +\lambda^2\tau f_tE_3 \\
  \overline{\nabla}_{\varphi_t}\varphi_t &=& 2\tau\lambda f_tE_1+(\lambda f_{tt}-\lambda^2f_t^2)E_2
\end{eqnarray*}
and,
\begin{eqnarray*}
 b_{11} &=& \lambda f_{xx}+\lambda^2(1+4\tau^2)f^2_x+2\tau\lambda^3(1+4\tau^2)f_xf_t+\lambda^2(1+4\tau^2) \\
 b_{12} &=& \lambda f_{xt}-\tau\lambda f^2_x+2\tau\lambda^3\left(\dfrac{1}{2}+2\tau^2\right)f_t^2-\tau\lambda \\
 b_{22} &=& \lambda f_{tt}-2\tau\lambda f_xf_t-\lambda^2f_t^2(1+4\tau^2)
\end{eqnarray*}
and,
\begin{eqnarray*}
 g_{11} &=& \lambda^2[(1+4\tau^2)+f_x^2] \\
 g_{12} &=& \lambda^2f_xf_t-2\tau\lambda \\
 g_{22} &=& 1+\lambda^2f^2_t
\end{eqnarray*}
a straightforward computation gives the result.
\end{proof}

An interesting formula for the Laplacian is given in the next Lemma.
\begin{lema}\label{l1.1}
 Considering $H=1/2$ then:\\
The mean curvature equation is given by
\begin{equation*}
 1=\dfrac{f^2}{W^3}[(f^2+f_t^2)f_{xx}-2(f_xf_t-2\tau f)f_{xt}+((1+4\tau^2)+f_x^2)f_{tt}+f(1+f_x^2)+2\tau f_xf_t]
\end{equation*}
And
\begin{eqnarray*}
 \Delta_S f &=& \dfrac{f^2}{W}\left(1-\dfrac{f}{W}+\dfrac{ff_x^2+2\tau f_tf_x}{W}\right) \\
 \Delta_S\left(\dfrac{1}{f}\right) &=& \dfrac{W-f}{fW}+\dfrac{f^2_t+2\tau(ff_xf_t+2\tau f_t^2)}{W}
\end{eqnarray*}

\end{lema}

\begin{proof}
 The proof follows from a hard computation by considering
\begin{equation*}
 \Delta_S=\dfrac{1}{\sqrt{g}}\displaystyle\sum_{ij}\partial_{x_{i}}(\sqrt{g}g^{ij}\partial_{x_{j}})
\end{equation*}
 where $g$ is the determinant of the first fundamental form and $(g^{ij})=(g_{ij})^{-1}$.

Observe that:
\begin{eqnarray*}
  \Delta_sf &=& \dfrac{1}{\sqrt{g}W^3}[f^2[(f^2+f_t^2)f_{xx}+2(2\tau f-f_xf_t)f_{xt}+(f^2_x+(1+4\tau^2))f_{tt}]+ \\
    & & +(a^3+f^3f_x)f_x+(af_x-(1+4\tau^2)ff_t)f_t]
\end{eqnarray*}
where $a=ff_x+2\tau f_t$ and $W^2=f^2+f^2_t+(ff_x+2\tau f_t)^2$
\end{proof}

In particular when $H=1/2$, we see from lemma \ref{l1.1}, that $f$ satisfies:
\begin{equation}\label{e1}
    (f^2+f^2_t)f_{xx}-2(f_xf_t-2\tau f)f_{xt}+(f^2_x+(1+4\tau^2))f_{tt}= -f(1+f^2_x)-2\tau f_xf_t+\dfrac{W^3}{f^2}
\end{equation}

Observe that we can rewriting equation (\ref{e1}) in the form:
\begin{equation*}
    (f^2+f^2_t)f_{xx}-2(f_xf_t-2\tau f)f_{xt}+(f^2_x+(1+4\tau^2))f_{tt}+
\end{equation*}
\begin{equation}\label{e2}
    +\left[2\tau f_t-f^2f_x\left(\dfrac{W^2}{f^2(W+f)}+\dfrac{1}{f}\right)\right]f_x -\left[\dfrac{W^2}{f^2(W+f)}+\dfrac{1}{f}\right](4\tau f_x+(1+4\tau^2)f_t)f_t=0
\end{equation}
Setting $a(x,t)=f+2\tau f_xf_t-\dfrac{W^3}{f^2}$. The equation (\ref{e2}) hold, by observing that:
\begin{eqnarray*}
  a(x,t) &=& -\left[\dfrac{W^3-f^3-2\tau f^2f_xf_t}{f^2}\right] \\
    &=& -\left[\dfrac{W^3-fW^2+fW^2-f^3-2\tau f^2f_xf_t}{f^2}\right] \\
    &=& -\left[\dfrac{W^2(W-f)}{f^2}+\dfrac{W^2-f^2-2\tau ff_xf_t}{f}\right] \\
    &=& -\left[\dfrac{W^2(W^2-f^2)}{f^2(W+f)}+\dfrac{W^2-f^2-2\tau ff_xf_t}{f}\right] \\
    &=& -\left[\dfrac{W^2}{f^2(W+f)}+\dfrac{1}{f}\right](W^2-f^2)+2\tau ff_xf_t \\
    &=& -\left[\dfrac{W^2}{f^2(W+f)}+\dfrac{1}{f}\right]((1+4\tau^2)f^2_t+f^2f^2_x+4\tau ff_xf_t)+2\tau ff_xf_t \\
    &=&  \left[2\tau f_t-\left(\dfrac{W^2}{f^2(W+f)}\right)f^2f_x\right]f_x-[(1+4\tau^2)f_t+4\tau ff_x] \left[\dfrac{W^2}{f^2(W+f)}+\dfrac{1}{f}\right]f_t
\end{eqnarray*}

\section{The First Theorem}

Before proving the first theorem, we prove the theorem \ref{T1}.

\begin{teorema}\label{T1} Let $U$ be the annulus $U=B_{R_2}\setminus B_{R_1}$ with $R_2\geq 2R_1$. Then for $\epsilon>0$ sufficiently small (depending only on $R_1$), there exist constant mean curvature $H=1/2$ horizontal graphs $f^+$ and $f^-$, satisfying equation (\ref{e1}) in $U$ with Dirichlet boundary data $f^{\pm}=1\pm\epsilon$ on $\partial B_{R_1}$, $f^{\pm}=1$ on $B_{R_2}$. Moreover $f^{\pm}$ tends to $1\pm\epsilon$ uniformly on compact subsets as $R_2$ tends to $\infty$.
\end{teorema}

\begin{proof} Let $U=B_{R_2}\setminus B_{R_1}$ be an annulus with $R_2\geq 4R_1$ and fix $h=1\pm\frac{\epsilon}{\log(R_2/R_1)}\log(\frac{R_2}{r})$, where $r^2=x^2+t^2$.


We expect the solution $f$ to be close to $h$, so we define the weighted $C^{2,\alpha}$ norm:
\begin{equation*}
    |v|^{*}_{2,\alpha;U}=\sup_{X}\{|v(X)|+r(X)|Dv(X)|+r^2(X)|D^2(X)|+r^{2+\alpha}_X[D^2v]_{\alpha;X}\}
\end{equation*}
where $X=(x,t)$ and $[D^2v]_{\alpha;X}$ is the $H\ddot{o}lder$ coefficient of $D^2v$ at $X$.
\begin{defi} We say $f$ is an admissible solution of (\ref{e1}) if $f\in \mathcal{A}_{\epsilon}$, where:
\begin{equation*}
    \mathcal{A}_{\epsilon}=\{f\in C^{2,\alpha}(U), f=h\hspace{0.2cm} on \hspace{0.2cm} U: |f-h|^{*}_{2,\alpha;U}\leq\sqrt{\epsilon}\}
\end{equation*}
\end{defi}

We note that $\mathcal{A}_{\epsilon}$ is convex and compact subset of the Banach space $\mathfrak{B}=C^{2,\beta}(U)$, $\beta<\alpha$. We will reformulate our existence problem as a fixed point of a continuous operator $T:\mathcal{A}_{\epsilon}\longrightarrow\mathcal{A}_{\epsilon}$ by rewriting  the equation (\ref{e2}) in the form:

\begin{equation*}
    (f^2+f^2_t)f_{xx}-2(f_xf_t-2\tau f)f_{xt}+(f^2_x+(1+4\tau^2))f_{tt}+
\end{equation*}
\begin{equation*}
    +\left[2\tau f_t-f^2f_x\left(\dfrac{W^2}{f^2(W+f)}+\dfrac{1}{f}\right)\right]f_x -\left[\dfrac{W^2}{f^2(W+f)}+\dfrac{1}{f}\right](4\tau f_x+(1+4\tau^2)f_t)f_t=0
\end{equation*}
\begin{obs} Note that this least equation implies that any solution $f^{\pm}$ solving the Dirichlet problem of Theorem \ref{T1} satisfies $1-\epsilon\leq f^-\leq1$ and $1\leq f^+\leq 1+\epsilon$ on $U$.
\end{obs}

We now define the operator $w=Tf$ as the solution of the linear Dirichlet problem

$\left\{
  \begin{array}{ll}
    L_fw:=aw_{xx}+2bw_{xt}+cw_{tt}+dw_x+ew_t=0, & \hbox{in U;} \\
       \hspace{0.65cm} w=h, & \hbox{on $\partial$U.}
  \end{array}
\right. $

where:
\begin{eqnarray*}
  a &=& f^2+f^2_x \\
  b &=& 2\tau f-f_xf_t \\
  c &=& f^2_x+(1+4\tau^2) \\
  d &=& ff_x+2\tau f_t-\left[\dfrac{W^2}{f^2(W+f)}+\dfrac{1}{f}\right]f^2f_x \\
  e &=& -[4\tau ff_x+(1+4\tau^2)f_t]\left[\dfrac{W^2}{f^2(W+f)}+\dfrac{1}{f}\right]
\end{eqnarray*}
\begin{defi} $D\subset\overline{U}$ is of scale $R$, if $X\in D$, then $c_1R\leq|X|\leq c_2R$ for uniform constant $c_1,c_2$.
\end{defi}

Note that for $f\in\mathcal{A}_{\epsilon}$, if $L_fu=F$ in $D\subset\overline{U}$ where $D$ is of scale $R$, then $\widetilde{u}=u(RX)$ satisfies:
\begin{equation}\label{e3}
    \widetilde{L}\widetilde{u}=\widetilde{a}\widetilde{u}_{xx}+2\widetilde{b}\widetilde{u}_{xt} +\widetilde{c}\widetilde{u}_{tt} +R\widetilde{d}\widetilde{u}_x+R\widetilde{e}\widetilde{u}_t=R^2\widetilde{F} \hspace{0.3cm} in\hspace{0.2cm} \overline{D}
\end{equation}
where $\widetilde{D}$ is of scale 1 and $\widetilde{a}(X)=a(RX)$, $\widetilde{b}(X)=b(RX)$, etc. Observe that, in this case $\widetilde{u}$ is defined on $\widetilde{D}=\{RX;c_1\leq|X|\leq c_2\}$.

Hence for $\epsilon$ sufficiently small, $\overline{L}$  is uniform close to $\Delta$ with $H\ddot{o}lder$ continuous coefficients.
\begin{prop} Let $w=Tf$ for $f\in\mathcal{A}_{\epsilon}$. Then for $\epsilon$ sufficiently small, $w\in\mathcal{A}_{\epsilon}$.
\end{prop}
\begin{proof} Set $u=w-h$, then:
\begin{equation}\label{e4}
    L_fu=[(1-f^2-f^2_t)h_{xx}+2f_xf_th_{xt}-f^2_xh_{tt}-dh_x-eh_t]:=F.
\end{equation}
By the maximum principle \cite[Theorem 3.1(pag. 32)]{GT}, $1\leq w\leq1+\epsilon$ (or $1-\epsilon\leq w\leq1$) so $|u|\leq\sqrt{\epsilon}$.

To see this, observe that for example $h^+=1+\dfrac{\epsilon}{\log(R_2/R_1)}\log(\frac{R_2}{r})$, so:
\begin{equation*}
    1\leq w \leq1+\epsilon \Leftrightarrow 1-h^+\leq u=w-h^+\leq 1+\epsilon+h^+
\end{equation*}
which is equivalent to:
\begin{equation*}
    -\epsilon\leq u \leq \epsilon
\end{equation*}
if and only if, $0\leq\frac{\log(R_2/r)}{\log(R_2/R_1)}\leq1$, which is true when $R_1\leq r\leq R_2$.
\vspace{0.5cm}

Since $R_2\geq4R_1$, there is a positive number $m_0$, such that, $R_2=m_0R_1$. Without loss of generality, we can suppose that $m_0$ is not a rational number. We now write $U=U_1\cup U_2\cup U_3$ where
\begin{eqnarray*}
 U_1 &=& \{X; R_1\leq|X|\leq\frac{m_0+2}{3} R_1\} \\
 U_2 &=& \{X; \frac{m_0+2}{3}R_1\leq|X|\leq2\frac{m_0+1}{3} R_1\} \\
 U_3 &=& \{X; 2\frac{m_0+1}{3}R_1\leq|X|\leq R_2\}
\end{eqnarray*}
Thus, each domain is of scale $R=R_1$ and we can apply Schauder interior or boundary estimates to $\widetilde{L}\widetilde{u}=R^2\widetilde{F}$ in $\widetilde{D}$ to obtain: (see \cite[Theorem 6.6(page 98)]{GT})
\begin{equation*}
    \|\widetilde{u}\|_{2,\alpha;\widetilde{D}}\leq C(\|\widetilde{u}\|_{0;\widetilde{D}}+\|R^2\widetilde{F}\|_{0,\alpha;\widetilde{D}}).
\end{equation*}
Observe that $|\widetilde{u}|\leq\epsilon$ implies $\parallel\widetilde{u}\parallel_{0;\widetilde{D}}\leq\epsilon$. From equation (\ref{e4}) follows $\|\widetilde{F}\|_{0,\alpha;\widetilde{D}}\leq C\epsilon^{\frac{3}{2}}$. This implies:
\begin{equation}\label{e5}
    \|\widetilde{u}\|_{2,\alpha;\widetilde{D}}\leq C(\|\widetilde{u}\|_{0;\widetilde{D}}+\|R^2\widetilde{F}\|_{0,\alpha;\widetilde{D}})\leq C\epsilon.
\end{equation}
This implies, (see \cite[Equation 4.17(page 61)]{GT})
\begin{equation*}
    \|\widetilde{u}\|^{*}_{2,\alpha;\widetilde{D}}\leq C\epsilon.
\end{equation*}

Undoing the scaling gives
\begin{equation*}
    \|u\|^{*}_{2,\alpha;\overline{D}}\leq C\epsilon.
\end{equation*}
Since $u=w-h$, it follows that for $\epsilon$ small enough, $w\in\mathcal{A}_{\epsilon}$ and the proposition is proved.

\end{proof}
We are now in a position to apply the Schauder fixed point theorem to our operator $w=Tf$ to find a solution $f^{\pm}\in\mathcal{A}_{\epsilon}$ to (\ref{e2}) which is equivalent to our original equation (\ref{e1}).

Now, we show that $f^{\pm}$ tends to $1\pm\epsilon$ uniformly on compact subsets as $R_2$ tends to $\infty$. To see this, observe that:
\begin{equation*}
    h^{\pm}(r)=1\pm\frac{\epsilon}{log(R_2/R_1)}\log(R_2/r)
\end{equation*}
tends to $1\pm\epsilon$ uniformly on compact subsets as $R_2$ goes to $\infty$. In fact, without loss of generality, we can suppose that $R_2=nR_1$, where $n$ is a positive integer (sufficiently large). Thus, we consider,
\begin{equation*}
    h_n(r)=\frac{\epsilon}{log(R_2/R_1)}\log(R_2/r)=\frac{\epsilon}{log(n)}\log(nR_1/r)=1+\frac{\log(R_1/r)}{\log(n)}
\end{equation*}
where $R_1\leq r$. Consider a compact domain $D$, then there exist constants $M_1,M_2>0$ such that,
\begin{equation*}
    -M_1\leq\log(R_1/r)\leq-M_2
\end{equation*}
for all $r\in[R_1,R_2]$ (closed interval). This implies,
\begin{equation*}
    \frac{-M_1}{\log(n)}\leq\frac{\log(R_1/r)}{\log(n)}\leq\frac{-M_2}{\log(n)}
\end{equation*}
Since, the function $\frac{1}{\log(n)}$ tends to $0$ uniformly on compact subsets, we conclude that, $h^{\pm}$ tends to $1\pm\epsilon$ uniformly on compact sets as $R_2$ tends to $\infty$. As $f^{\pm}$ is $C^{2,\alpha}$ close to $h^{\pm}$, we conclude the second affirmation of the theorem.

\end{proof}

\begin{teorema} Let $\Sigma$ be a properly embedded constant mean curvature $H=1/2$ surface in $\widetilde{PSL}_2(\R,\tau)$. Suppose $\Sigma$ is asymptotic to a horocylinder $C$, and one side of $C$. If the mean curvature vector of $\Sigma$ has the same direction as that of $C$ at points of $\Sigma$ converging to $C$, then $\Sigma$ is equal to $C$ (or a subset of $C$ if $\partial\Sigma\neq\phi$).
\end{teorema}
\begin{proof} After an isometry, we can assume that, there is a sequence of points $p_i=(x_i,y_i,t_i)\in\Sigma$ with $y_i\rightarrow1$ and $\Sigma$ is asymptotic to $C(1)$ when $t$ goes to $+\infty$ (here $C(\eta)$ denote the $1/2$ horocylinder $y=\eta$). In this case $\langle\overrightarrow{H},\partial_y\rangle>0$. Thus, either $\Sigma$ is contained in the set $\{y<1\}$ or $\Sigma$ is contained in the set $\{y>1\}$. First, we suppose that $\Sigma$ is contained in the set $\{y<1\}$.

For $\epsilon>0$ we consider the slab $S^-$ bounded by $C(1-\epsilon)$ and $C(1)$, then by the maximum principle $\Sigma^-=\Sigma\cap S^-$ has a non compact component with boundary $\partial\Sigma\subset C(1-\epsilon)$.

Let $D(\eta,R)$ denote the disk in $C(\eta)$ defined by $D(\eta,R=\{(x,\eta,t);x^2+t^2\leq R^2\}$. By considering vertical translation, we can find a disk $D(1,3R_1)$ such that:
\begin{equation*}
    D(1,3R_1)\times[1-\epsilon,1]\cap\Sigma^-=\phi.
\end{equation*}
By Theorem \ref{T1}, for each $R\geq2R_1$, there exist a horizontal graph $f^-_R$ defined on the annulus $U=B_{R_1}\ B_{R_2}$, this horizontal graph converge to $C(1-\epsilon)$, when $R$ goes to $+\infty$.

Now, consider $R$ large, such that the graph of $f^-_R$ (which we denote by $\Gamma^-$), satisfies $\Sigma^-\cap\Gamma^-\neq\phi$. by considering translations along the geodesic $\{x=0,t=0\}$ together with vertical translation, the translated surface of $\Gamma^-$ does not touch $\Sigma^-$, that is, there is a translated surface of $\Gamma^-$ (which we denote by $\Gamma^-_1$) such that $\Gamma^-_1$ and $\Sigma^-$ has an interior contact point. Since the mean curvature vectors are pointing up, his violates the maximum principle and $\Sigma^-$ cannot exist.

In the second case, we redo exactly the same argument exchanging the roles of $C(1+\epsilon)$ and $C(1-\epsilon)$.

\end{proof}

\section{H Sections in $\widetilde{PSL}_2(\R,\tau)$}

\indent

We give the definition of Killing Submersion which is due to $H$. $Rosenberg$, $R$. $Souam$, and $E$. $Toubiana$ see \cite{RRT}.

Consider a Riemannian 3-manifolds $(M^{3},g)$ which fibers over a Riemannian surface $(M^{2},h)$, where $g$ and $h$ denote the Riemannian metrics respectively

\begin{defi}\label{a1} A Riemannian submersion $\pi:(M^{3},g)\longrightarrow (M^{2},h)$ such that:
\begin{enumerate}
  \item each fiber is a complete geodesic,
  \item the fibers of the fibration are the integral curves of a unit Killing vector field $\xi$ on $M^{3}$.
\end{enumerate}
will be called a Killing submersion.
\end{defi}
\begin{defi}\label{a2} Let $\pi:(M^{3},g)\longrightarrow(M^{2},h)$ be a Killing submersion.
\begin{enumerate}
  \item Let $\Omega\subset M^{2}$ be a domain. An $H$-$section$ over $\Omega$ is an $H$-$surface$ which is the image of a section.
  \item Let $\gamma\subset M^{2}$ be a smooth curve with geodesic curvature 2H. Observe that the surface $\pi^{-1}(\gamma)\subset M^{3}$ has mean curvature H. We call such a surface a $vertical$ $H-cylinder$.
\end{enumerate}
\end{defi}

Observe that, $\pi:\widetilde{PSL}_2(\R,\tau)\longrightarrow \mathbb{H}^{2}$ given by $\pi(x,y,t)=(x,y)$ is a Killing submersion. \\
Since, $\pi:\widetilde{PSL}_2(\R,\tau)\longrightarrow \mathbb{H}^{2}$ is a Killing submersion, we can consider graphs in $\widetilde{PSL}_2(\R,\tau)$.
\begin{defi} A graph in $\widetilde{PSL}_2(\R,\tau)$ over a domain $\Omega$ of $\mathbb{H}^2$ is the image of a section $s_0:\Omega\subset\mathbb{H}^2\longrightarrow\widetilde{PSL}_2(\R,\tau)$.
\end{defi}
Given a domain $\Omega\subset\mathbb{H}^2$ we also denote by $\Omega$ its lift to $\mathbb{H}^2\times0$, with this identification we have that the graph (vertical graph) $\Sigma(u)$ of $u\in(C^0(\partial\Omega)\cap C^{\infty}(\Omega))$ is given by:
\begin{equation*}
    \Sigma(u)=\{(x,y,u(x,y))\in\widetilde{PSL}_2(\R,\tau);(x,y)\in\Omega\}
\end{equation*}

\section{The Second Theorem}

\indent

In this section our second main result concerns complete $H=1/2$ surfaces in $\widetilde{PSL}_2(\R,\tau)$ transverse to the vertical Killing field $E_3=\partial_t$. We prove such surfaces are entire graphs.

\begin{teorema} Let $\Sigma$ be a complete immersed surface in $\widetilde{PSL}_2(\R,\tau)$ of constant mean curvature $H=1/2$. If $\Sigma$ is transverse to $E_3$ then $\Sigma$ is an entire vertical graph over $\mathbb{H}^2$.
\end{teorema}
\begin{proof} The proof is the same as this one to $\mathbb{H}^2\times\R$ (see \cite{HRS}). By completeness we give it. From now on we identify $\mathbb{H}^2$ with its lift $\mathbb{H}^2\times\{0\}$\\
The mean curvature vector of $\Sigma$ never vanish so $\Sigma$ is orientable. Let $\nu$ be a unit vector field along $\Sigma$ in $\widetilde{PSL}_2(\R,\tau)$. The function $u=\langle\nu,E_3\rangle$ is a non-zero Jacobi function on $\Sigma$, so $\Sigma$ is strongly stable and thus has bounded curvature. We can assume $u>0$ and $\langle\nu,\overrightarrow{H}\rangle>0$

Here there is $\delta>0$ such that for each $p\in\Sigma$, $\Sigma$ is a graph (in exponential coordinates) over the disk $D_{\delta}\subset T_p\Sigma$ of radius $\delta$, centered at the origin of $T_p\Sigma$. This graph, denoted by $G(p)$, has bounded geometry. The $\delta$ is independent of $p$ and the bound on the geometry of $G(p)$ is uniform as well (see \cite{RRT}).


We denote by $F(p)$ the surface $G(p)$ translated to the origin $i\in\mathbb{H}^2\approx\mathbb{H}^2\times\{0\}$. (The translation that takes $p$ to $i$).


For $q\in\mathbb{H}^2$, we denote by $\Gamma_{\delta}(q)$ a horizontal horocycle arc of length $2\delta$, centered at $q$.\newline
\\

$\mathbf{Claim 1}$: Let $p_n\in\Sigma$, satisfy $u(p_n)\rightarrow0$ as $n\rightarrow\infty$ ($T_{p_n}(\Sigma)$ are becoming vertical). There is a subsequence of $p_n$ (which we also denote by $\{p_n\}$) such that $F(p_n)$ converges to $\Gamma_{\delta}(i)\times[-\delta,\delta]$, for some horocycle $\Gamma_{\delta}(i)$. The convergence is in the $C^2$-topology.
\\

Proof of Claim 1. Choose a subsequence $p_n$ so that the oriented tangent planes $T_i(F(p_n))$ converge to a vertical plane $P$. Let $\Gamma_{\delta}(i)$ be the horocycle arc through $i\in\mathbb{H}^2\approx\mathbb{H}^2\times0$ whose curvature vector has the same direction as the curvature vector of the (limit) curvature vectors of $F(p_n)$.

Since the $F(p_n)$ have bounded geometry and they are graphs over $D_{\delta}(p_n)\subset T_{p_n}(F(p_n))$, the surfaces $F(p_n)$ are bounded horizontal graphs over $\Gamma_{\delta}(i)\times[-\delta,\delta]$ for $n$ large. Thus a subsequence of these graphs converges to an $H=1/2$ surface $F$; $F$ is tangent to $\Gamma_{\delta}(i)\times[-\delta,\delta]$ at $i$ and a horizontal graphs over this. It suffices to show $F=\Gamma_{\delta}(i)\times[-\delta,\delta]$.

Were this not the case, then the intersection near $i$, of $F$ and $\Gamma_{\delta}(i)\times[-\delta,\delta]$ would consist of $m$ smooth curves passing through $i$, $m\geq2$, meeting transversally at $i$. In a neighborhood of $i$, theses curves separate $F$ into $2m$ components. Adjacents components lie on opposite sides of $\Gamma_{\delta}(i)\times[-\delta,\delta]$.

Hence in a neighborhood of $i\in F$, the mean curvature vector of $F$ alternates from pointing up in $\widetilde{PSL}_2(\R,\tau)$ to pointing down (or vice-versa), as one goes from one component to the other. But $F(p_n)$ converges to $F$ in the $C^2$-topology, so $F(p_n)$, $n$ large, would also have points where the mean curvature vector point up and down in $\widetilde{PSL}_2(\R,\tau)$. This contradicts that $F(p_n)$ is transverse to $E_3$, and claim 1 is proved. Notice that we have proved that whenever $F(p_n)$ converges to a local surface $F$, $F$ is necessarily some $\Gamma_{\delta}(i)\times[-\delta,\delta]$. This prove Claim 1.
\\

Now let $p\in\Sigma$ and assume $\Sigma$ in a neighborhood of $p$ is an $H$ section, that is, in a neighborhood of $p$, $\Sigma$ is a vertical graph of a function $f$ defined on $B_R$, $B_R$ the open ball of radius $R$ of $\mathbb{H}^2$, centered at $i\in\mathbb{H}^2$. Denote by $S(R)$ the graph of $f$ over $B_R$. If $\Sigma$ is not an entire graph then we let $R$ be the largest such $R$ so that f exist. Since $\Sigma$ has constant mean curvature, $f$ has bounded gradient on relatively compact subsets of $B_R$, see \cite[Theorem 3.5]{RRT}.

Let $q\in\partial B_R$ be such that $f$ does not extend to any neighborhood of $q$ (to an $H=1/2$ graph).
\\

$\mathbf{Claim 2}$: For any sequence $q_n\in B_R$, converging to $q$, the tangent planes $T_{p_n}(S(R))$, $p_n=(q_n,f(q_n))$, converge to a vertical plane $P$. $P$ is tangent to $\partial B_R$ at $q$ (after vertical translation to height zero in $\widetilde{PSL}_2(\R,\tau)$).
\\

Proof of Claim 2. Let $F(n)$ denote the image of $G(p_n)$ under the vertical translation taking $p_n$ to $q_n$. Observe first, that $T_{q_n}(F(n))$ converges to the vertical , for any subsequence of the $q_n$. Otherwise the graph of bounded geometry $G(p_n)$, would extend to a vertical graph beyond $q$, for $q_n$ close enough to $q$, hence $f$ would extend; a contradiction.

Now we can prove $T_{q_n}(F_n)$ converges to the vertical plane $P$ passing  through $q$ and tangent to $\partial B_R$ at $q$. Suppose some subsequence $q_n$ satisfies $T_{q_n}(F_n)$ converges to a vertical plane $Q$, $Q\neq P$, $q\in Q$. By Claim 1, the $F_n$ converge in the $C^2$-topology, to $\Gamma_{\delta}\times[-\delta,\delta]$, where $\Gamma_{\delta}(q)$ is a horocycle arc centered at $q$. Since $Q\neq P$, and $\Gamma_{\delta}$ is tangent to $Q$ at $q$, there are points of $\Gamma_{\delta}$ in $B_R$. Such a point is the limit of points on $F_n$. Then the gradient at $f$ at these points of $F_n$ diverges, which contradicts interior gradient estimates of $f$. This proves Claim 2.
\\

Now applying Claim 1 and Claim 2, we know that for any sequence $q_n\in B_R$ converging to $q$, the $F(q_n)$ converge to $\Gamma_{\delta}(q)\times[-\delta,\delta]$.
\\

$\mathbf{Claim 3}$: For any $q_n\rightarrow q$, $q_n\in B_R$, we have $f(q_n)\rightarrow+\infty$ or $f(q_n)\rightarrow-\infty$.
\\

Proof of Claim 3. Let $\gamma$ be a compact horizontal geodesic of length $\epsilon$ starting at $q$, entering $B_R$ at $q$, and orthogonal to $\partial B_R$ at $q$. Let $C$ be the graph of $f$ over $\gamma$. Notice that $C$ has no horizontal tangents at points near $q$ since the tangent planes of $S(R)$ are converging to $P$. So assume $f$ is increasing along $\gamma$ as one converges to $q$. If $f$ were bounded above, then $C$ would have a finite limit point $(q,c)\in\Sigma$. But then $\Sigma$ has a vertical tangent plane at $(q,c)$; a contradiction. This prove Claim 3.
\\

Now choose $q_n\in\gamma$, $q_n\rightarrow q$, and $F(q_n)$ converges to $\Gamma_{\delta}(q)\times[-\delta,\delta]$. Let $\Gamma$ be the horocycle containing $\Gamma_{\delta}(q)$, and parameterize $\Gamma$ by arc length; denote $q(s)\in\Gamma$ the point at distance $s$ on $\Gamma$ from $q=q(0)$, $-\infty<s<+\infty$. Denote by $\gamma(s)$ a horizontal geodesic arc orthogonal to $\Gamma$ at $q(s)$, $q(s)$ the mid-point of $\gamma(s)$. Assume the length of each $\gamma(s)$ is $2\epsilon$ and $\cup_{s\in\R}\gamma(s)=N_{\epsilon}(\Gamma)$ is the $\epsilon$-tubular neighborhood of $\Gamma$.

Let $\gamma^+(s)$ be the part of $\gamma(s)$ on the mean convex side of $\Gamma$; so $\gamma=\gamma^+(0)$. More precisely, the mean curvature vector of $\Sigma$ points up in $\widetilde{PSL}_2(\R,\tau)$, and $f\rightarrow+\infty$ as one approaches $q$ along $\gamma$, so $\Gamma$ is convex towards $B_R$.
\\

$\mathbf{Claim 4}$: For $n$ large, each $F(q_n)$ is disjoint from $\Gamma\times\R$. Also, for $|s|\leq\delta$, $F(q_n)\cap(\gamma^+(s)\times\R)$ is a vertical graph over an interval of $\gamma^+(s)$.
\\

Proof of Claim 4. Choose $n_0$ such that for $n\geq n_0$, $C_n(s)=F(q_n)\cap(\gamma(s)\times\R)$ is one connected curve of transverse intersection, for each $s\in[-\delta,\delta]$. Since the $F(q_n)$ are $C^2$-close to $\Gamma_{\delta}(q)\times[-\delta,\delta]$, $C_n(s)$ has no horizontal or vertical tangents and is a graph over an interval in $\gamma(s)$.

We now show this interval is in $\gamma^+(s)-q(s)$. Suppose not, so $C_n(s)$ goes beyond $\Gamma\times\R$ on the concave side. Recall that $C=\gamma\cap P^{\perp}$ is the graph of $f$ and $f\rightarrow+\infty$ as one goes up on $C$. We have $p_n=(q_n,f(q_n))$. Fix $n\geq n_0$ and choose new points $q_k$, $k\geq n$, so that $f(q_{k+1})-f(q_k)=\delta$; clearly $q_k\rightarrow q$ as $k\rightarrow+\infty$. Lift each $C_k(s)$ to $G(p_k)$ by the vertical translation of $F(q_k)$ by $f(q_k)$. By construction, $C_{k+1}(s)$ is the analytic continuation of $C_k(s)$ in $\Sigma\cap(\gamma(s)\times\R)$, for each $s\in[-\delta,\delta]$ and for all $k\geq n+1$. The curve $C(s)=\cup_{k\geq n}C_k(s)$ is a vertical graph over an interval in $\gamma(s)$. It has points on the concave side of $\Gamma\times\R$ for some $s_0\in[-\delta,\delta]$. For $s=0$, $C(0)=C$ stays on the convex side of $\Gamma\times\R$. So for some $s_1$, $0<s_1\leq s_0$, $C(s_1)$ has a point on $\Gamma\times\R$ and also inside the concave side of $\Gamma\times\R$.

But the $F(q_k)$ converge uniformly to $\Gamma_{\delta}(q)\times[-\delta,\delta]$ as $k\rightarrow\infty$, so the curve $C(s_1)$ converges to $q(s_1)\times\R$ as the height goes to $\infty$, This obliges $C(s_1)$ to have a vertical tangent on the concave side of $\Gamma\times\R$, a contradiction. This prove Claim 4.
\\

Now we choose $\epsilon_1<\epsilon$ (which we call $\epsilon$ as well) so that $\cup_{s\in[-\delta,\delta]}C(s)$ is a vertical graph of a function $g$ on $\cup_{s\in[-\delta,\delta]}(\gamma^+(s)-q(s))$, (the $\gamma^+(s)$ now have length $\epsilon_1$); $g$ is an extension of $f$.

The graph of $g$ on each $\gamma^+\times\R$ is the curve $C(s)$, and the graph of $g$ converges to $\Gamma_{\delta}(q)\times\R$ as the height goes to infinity.

Now we begin this process again replacing $C$ by the curves $C(\delta)$ and then $C(-\delta)$. Analytic continuation yields an extension $h$ of $g$ to a domain $\Omega$ contained in the open $\epsilon$-tubular neighborhood of $\Gamma\times\R$, on the convex side of $\Gamma$. $\Omega$ is an open neighborhood of $\Gamma$ in this mean convex side. The graph $h\rightarrow\infty$ as one approaches $\Gamma$ in $\Omega$; it converges to $\Gamma\times\R$ as the height goes to infinity.
\\

$\mathbf{Claim 5}$: There is an $\epsilon>0$, such that $\Omega$ contains the $\epsilon$ tubular neighborhood of $\Gamma$ on the convex side.
\\

Proof of Claim 5. We know there is a domain $\Omega$ on the convex side of $\Gamma$; $\Omega$ is a neighborhood of $\Gamma$ on the convex side. Also the surface $\Sigma$ contains a graph over $\Omega$, composed of curves $C(q)$, $q\in\Gamma$, where each curve $C(q)$ is a graph over an interval $\gamma^+(q)$, $\gamma^+(q)$ orthogonal to $\gamma$ at $q$. Also $C(q)$ is a strictly monotone increasing graph with no horizontal tangents and $C(q)$ converges to $\{q\}\times\R^+$, as one goes up to $+\infty$.


The graph over $\Omega$ is converging uniformly to $\Gamma\times\R^+$ as one goes up.

Now suppose that for some $q\in\Gamma$, $\gamma^+(q)$ is of length less than $\epsilon$. Then $C(q)$ diverges to $-\infty$ as one approaches the end-point $\widetilde{q}$ of $\gamma^+(q)$, $\widetilde{q}\neq q$.


The previous discussion where we showed the graph over $\Omega$ exists and converges to $\Gamma\times\R^+$, now applies to show that there is a horocycle $\widetilde{\Gamma}$ passing through $\widetilde{q}$, $C(q)$ converges to $\{\widetilde{q}\}\times\R^{-}$ as one tends to $\widetilde{q}$ on $\gamma^+(q)$. Also a $\delta$-neighborhood of $C(q)$ in $\Sigma$, converges uniformly to $\widetilde{\gamma}_{\delta}(\widetilde{q})\times\R^-$, as one goes down to $-\infty$. We know this $\delta$-neighborhood of $C(q)$ in $\Sigma$, converges uniformly to $\Gamma_{\delta}(q)\times\R^+$, as one goes up to $+\infty$.

For each $q(s)\in\Gamma$, a distance $s$ from $q$ on $\Gamma$, $|s|\leq\delta$, the curve $C(q(s))$ converges uniformly to some $\{\widetilde{q}(\widetilde{s})\times\R^-\}$, as one goes down to $-\infty$. By analytic continuation of the $\delta$-neighborhoods, one continues this process alon $\gamma$.

If $\Gamma\cup\widetilde{\Gamma}=\phi$, then the process continues along all $\Gamma$ and $\Omega$ is the region bounded by $\Gamma\cup\widetilde{\Gamma}$. This suffices to prove Claim 5 since each $\gamma^+(q)$, $q\in\Gamma$, has the same length.

So we can assume $\Gamma\cup\widetilde{\Gamma}=\{p\}$. Consider the curves $C(q(s))$, as $q(s)$ goes from $q$ to $p$ along $\Gamma$. They are graphs that become vertical both at $+\infty$ and $-\infty$. Hence the graphs $C(q(s))$ become vertical at every point as $q(s)\rightarrow p$.


Consider the point of $C(q(s))$ at height $0$ in $\widetilde{PSL}_2(\R,\tau)$. As $q(s)\rightarrow p$, these converge to a point of $\Sigma$ and the tangent plane of $\Sigma$ is vertical at this point; a contradiction.

We remark that in this case $f(q_n)\rightarrow-\infty$ (see Claim 3), one works on the concave side of the horocycle $\Gamma(q)$ and Claims 4 and 5 show there is an $\epsilon>0$ and a graph $G\subset\Sigma$ over the domain $\Omega(\epsilon)$ between $\Gamma(\epsilon)$ (the equidistant horocycle to $\Gamma$ on the concave side of $\Gamma$) and $\Gamma$. The graph $G$ converges uniformly to $\Gamma\times\R$ as one approaches $\Gamma$ in $\Omega(\epsilon)$.

To complete the proof of the theorem we apply the half-space theorem to show no such $H=1/2$ graph exist. Strictly speaking we can not apply the half-space theorem directly since the $G$ does not have a compact boundary. But this graph is proper in the tubular neighborhood of the horocylinder, so the proof shows the graph can not exist.

\end{proof}


\begin{thebibliography}{100}
\bibitem{D} Daniel, B. Isometric immersions into 3-dimensional homogeneous manifolds. Comment. Math. Helv. 82 (2007), no. 1, 87-131.
\bibitem {RRT} Rosenberg, H.; Souam, R.; Toubiana, E.. General Curvature Estimates for Stable H-Surfaces in 3-Manifolds. To appear in J. of Diff. Geom.
\bibitem {HRS} Hauswirth, L.; Rosenberg, H.; Spruck, J.. On Complete Mean Curvature 1/2 surfaces in $\mathbb{H}^{2}\times\R$. Comm. Anal. Geom., 16 (5) (2009) 989-1005
\bibitem {GT} Gilbarg, D.; Trudinger, N. S.. Elliptic Partial Differenial Equations of Second Order. Second edition, 1997, Springer.
\bibitem {WT} Thurston, W. Three-Dimensional Geometry and Topology. Princeton, 1997.


\end{thebibliography}
\end{document}